\documentclass[conference]{IEEEtran}
\usepackage{caption}
\IEEEoverridecommandlockouts
\usepackage{cite}
\usepackage{amsmath,amssymb,amsfonts}
\usepackage{graphicx}
\usepackage{algorithm}
\usepackage{algpseudocode} 
\usepackage{textcomp}
\usepackage{xcolor}
\usepackage{booktabs}
\usepackage{subcaption}
\usepackage{adjustbox}
\usepackage{tikz}
\usepackage{pgfplots}
\pgfplotsset{compat=1.18}
\usepackage{amsmath}
\usepackage[american,siunitx]{circuitikz}
\usetikzlibrary{positioning}
\usepackage{mathptmx}
\usepackage{newtxtext}
\usepackage{url} 

\def\BibTeX{{\rm B\kern-.05em{\sc i\kern-.025em b}\kern-.08em
    T\kern-.1667em\lower.7ex\hbox{E}\kern-.125emX}}

\def\addlegendimage{\csname pgfplots@addlegendimage\endcsname}

\begin{document}

\newif\ifdraft
\drafttrue 

\ifdraft
        \newcommand{\diana}[1]{\textcolor{violet}{[[Diana: #1]]}}
        \newcommand{\nicolasc}[1]{\textcolor{brown}{[[Nicolas: #1]]}}
        \newcommand{\soummya}[1]{\textcolor{teal}{[[Soummya: #1]]}}
\else
        \newcommand{\diana}[1]{}
        \newcommand{\nicolasc}[1]{}
        \newcommand{\soummya}[1]{}
\fi
 
\title{Peer-to-Peer (P2P) Electricity Markets for Low Voltage Networks\\
\thanks{This work was partially funded by Funda\c{c}\~{a}o para a Ci\^{e}ncia e a Tecnologia under the scope of the Carnegie Mellon Portugal Program.}
}
 
\author{\IEEEauthorblockN{Diana Vieira Fernandes}
\IEEEauthorblockA{\textit{Dept. of Engineering \& Public Policy} \\
\textit{Carnegie Mellon University, USA}\\
\textit{IN+/LARSyS, Instituto Superior T\'ecnico} \\
\textit{Universidade de Lisboa, Portugal} \\
\texttt{dianaimf@andrew.cmu.edu}}
\and
\IEEEauthorblockN{Nicolas Christin}
\IEEEauthorblockA{\textit{Dept. of Engineering \& Public Policy} \\
\textit{Carnegie Mellon University, USA}\\
\texttt{nicolasc@andrew.cmu.edu}}
\and
\IEEEauthorblockN{Soummya Kar}
\IEEEauthorblockA{\textit{Dept. of Electrical \& Computer Engineering} \\
\textit{Carnegie Mellon University, USA}\\
\texttt{soummyak@andrew.cmu.edu}}
}

\maketitle
\begin{abstract}
We develop a clearance and settlement model for Peer-to-Peer
(P2P) energy trading in low-voltage networks. The model enables direct
transactions between parties within an open and distributed system and
integrates unused capacity while respecting network
constraints. We evaluate the model through simulations of different scenarios (normal operating conditions and extreme conditions) for 24-hour time blocks. Our simulations highlight the benefits of our model in a decentralized energy system, notably its
ability to deal with high-trade volumes.
\end{abstract}


\section{Introduction}
The increasing penetration of decentralized renewable generation sources, such as solar photovoltaics (PV), along with the digitization of the electric grid, opens the door for innovative market structures and mechanisms, namely, peer-to-peer (P2P) \cite{morstyn_using_2018, zhang_distributed_2019,
mengelkamp_designing_2018, guerrero_decentralized_2019,
azim_regulated_2020, ullah_peer--peer_2021, soto_peer--peer_2021,
capper_peer--peer_2022, guo_online_2021, sousa_peer--peer_2019}, besides ``energy communities'' and ``local energy markets''~\cite{capper_peer--peer_2022}. Low-voltage networks, which are prevalent in distribution systems, offer untapped potential for such markets. However, the lack of an established mechanism for secure and efficient trading and operations is a significant challenge and the implementation has yet to be fully realized.

Peer-to-peer energy trading could democratize the energy market (granting effective access rights, and participatory mechanisms), where each user can either be a consumer or producer, i.e. a ``prosumer,'' reduce energy losses (minimizing the need for long-distance electricity transportation \cite{zhang_distributed_2019}), and foster competition by enabling broader market participation.

Traditionally, the electricity sector relied on a vertically integrated utility model. The utility procured electricity to meet the total demand (aggregated load) of all consumers. The (centralized) Economic Dispatch (ED) model presumes an aggregated load (or demand), focusing solely on wholesale or large-scale generation. Its primary objective is to minimize the costs associated with electricity generation and distribution \cite{carpentier_contribution_1962,carpentier_optimisation_1963} given the (aggregated) power output (and generation cost) of each generator. As such, clearance and settlement procedures were primarily directed toward the Independent System operator/Transmission System Operator (ISO/TSO). Distribution System Operators (DSOs) were overlooked, as primarily handled the distribution of electricity, from retailers to consumers. The latter did not have generation capacity, and the ability to put in the market - excess electricity - is changing the traditional market structures and roles. However, electricity markets modeling still use a conventional, hierarchical approach in power system management, resulting in ``prosumers'' typically assuming the role of passive participants \cite{sousa_peer--peer_2019}, failing to capture the complexities of modern energy distribution systems and regulations mandating the nondiscriminatory access to the network (e. g. section 206 (b) of the Federal Power Act (FPA) \cite{noauthor_federal_1935} (wholesale markets) or art. 3º of the (EU) Directive 2019/944 on common rules for the internal electricity market \cite{european_commission_directive_2019}.

Peer-to-peer (P2P) electricity trading refers to the direct transaction
of electricity between two parties (producers and consumers) within
the electricity network without the intermediation of traditional
centralized or aggregating entities (i.e. retail companies). In essence,
P2P is a bilateral contract, i.e., a micro Power Purchase Agreement
(PPA), by which one party undertakes to place the contracted electricity 
on the grid and the
other to receive it, at prices and conditions
freely agreed (i.e. PJM (US) market operations \cite{pjm_pjm_2023,
pjm_pjm_2024} and Portugal \cite{erse_manual_2023}). In electricity
operations, whether engaging in bilateral contracts or organized
markets, all injected electricity is treated equally (as another
fungible commodity) for dispatch and clearance. The classical
distinction between wholesale and retail is also blurred. Thus, ED-based
models are not well-aligned with real-world scenarios for decentralized
decision-making and distributed networks, where energy distribution
and dispatch are driven by local conditions (local transactions),
bidirectional flows, network-specific constraints, and actual network
usage, particularly when multiple -- not aggregated -- small-scale
``prosumers'' can directly sell to others and use the network
simultaneously and actively.
Considering network operations, it is assumed that several contracting schemes coexist (market, bilateral transaction and other, as currently are already established in most operation regulations -- e.g. in the US, for the PJM area \cite{pjm_pjm_2023, pjm_pjm_2024} and Portugal \cite{erse_manual_2023}) where there is an integration for of several transactions (orders) into the network power flow management (independent of the underlying contracting arrangements). Several transactions occur concurrently on the same network for dispatch with different time horizons and maturities. 
This work addresses the transformation of traditional market structures by allowing consumers with generation capacity to offer excess electricity (that can be direct or into a pool), challenging prior assumptions about power flow modeling of P2P operations. Through this approach, this work elucidates the potential impacts of P2P trading on power distribution networks, highlighting the need for novel models to understand and incorporate the emerging electricity market and operation dynamics. To simulate these transactions in realistic or near-realistic conditions, a double auction process is used (for the matching process) of potential P2P transactions and mirrors operations within wholesale markets (managed by the TSO that performs power flow analyses to manage and build dispatch schedules, thereby integrating and coordinating various contracting schemes and temporal resolutions) to accurately model the nuanced and localized nature of P2P, thereby mapping the flow of electricity and contracts. Existing works have highlighted the potential benefits of peer-to-peer (P2P) energy but have not modeled these operations integrated into actual network operations, and typically do not include the physical constraints~\cite{capper_peer--peer_2022}. Most simulations assume a closed network, particularly within microgrids \cite{mengelkamp_designing_2018, zhang_distributed_2019} or within virtual power plants \cite{morstyn_using_2018}, where all users are part of the network or system set-up, focusing primarily on welfare maximization \cite{azim_regulated_2020} or aggregated cost minimization
\cite{guo_online_2021, ullah_peer--peer_2021, sousa_peer--peer_2019}, similar to the approach to the (centralized) ED problem. These models simplify the trading environment by assuming a controlled group of participants -- most focus on small energy markets with less than 10 participants \cite{capper_peer--peer_2022}. In reality, electricity networks are generally open and interconnected. Not all users participate in P2P trading (as not all consumers contract from the same retail company), and traditional producers and retailers coexist \cite{hug_consensus_2015}. 
The DSO’s role has expanded to include power flow coordination and control due to the introduction of generation capacities at the distribution level (e.g. net metering). This role is distinct from -- retail or wholesale -- generation or intermediate entities such as retail companies that operate the connection between markets. In some cases, the complete separation of these roles may be less pronounced or even bundled. 

\subsection{Contributions}
We formalize and develop a novel model for a peer-to-peer clearance and settlement problem designed for integrating peer-to-peer energy trading in low-voltage, incorporating network constraints with multiple and concurrent autonomous users and different time windows. To the best of our knowledge, this work addresses aspects not comprehensively covered in existing and prior works, which:

\begin{itemize}
    \item Fills a research gap by integrating distributed power generation of small, heterogeneous, and autonomous agents within Peer-to-Peer (P2P) trading mechanisms;
    \item Allows the integration of unused close capacity and excess surplus renewable energy;
    \item Facilitates the execution of viable P2P energy exchanges, alongside traditional electricity retailers, merging these resources with existing market and operations mechanisms, in an open manner, with multiple users and transactions;
    \item Promotes the evolution of novel market configurations -- agnostic to contracting mechanisms, and does not assume prior relationships between participants -- facilitating the engagement of ``prosumers'' in energy trading while ensuring power flow remains within permissible bounds;
    \item Power flow equations are explicitly integrated where network topology is considered (to reflect the network state from several concurrent transactions);
    \item Allows clearance and integration with existing market structures and outputs a pre-operational dispatch schedule (assuming a centralized DSO that receives all potential trades from independent participants);
    \item Incorporates network constraints and considers power flow in trade execution: unlike financial contracting, the trade occurs only if it can be cleared and physically settled (or can be dispatched and delivered); 
    \item Preliminary testing on various trading scenarios and networks demonstrate our model can facilitate and integrate P2P power flow within acceptable network constraints, even in extreme cases.
\end{itemize}

\section{Methods -- System Setup} \label{sec:methods}
\begin{figure*}[htpb]
\vspace{-2mm}
  \centering
  \includegraphics[width=\textwidth]{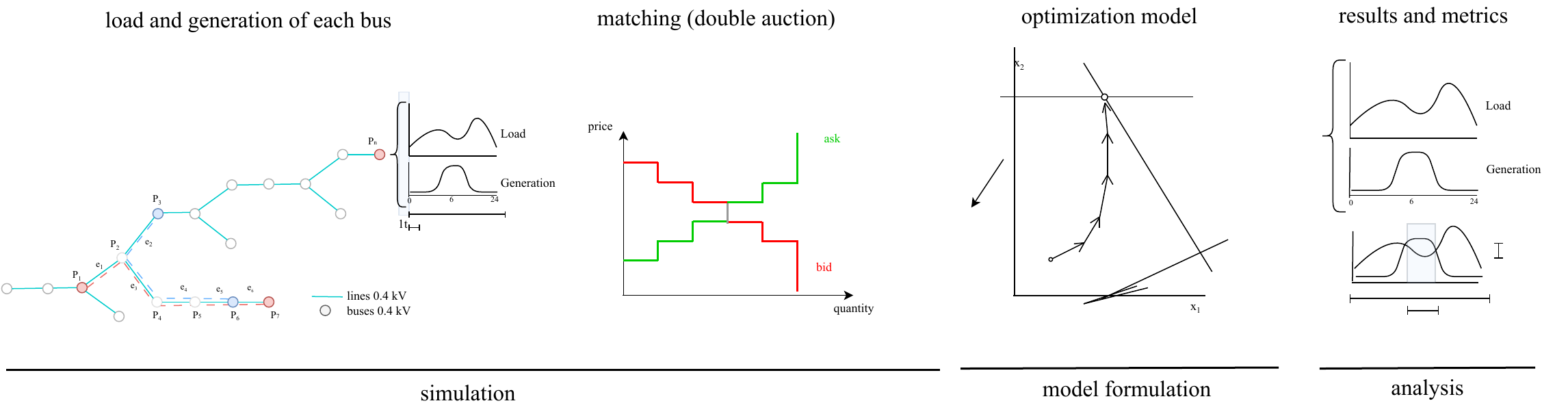}
  \caption{Methods and processes}
  \label{fig:methods}
\vspace{-1mm}
\end{figure*}

In the context of distributed energy systems, there is a paradigm shift in how energy is produced, distributed, and consumed. This section describes our approach to modeling and analyzing peer-to-peer (P2P) transactions -- as a two-stage process -- emphasizing the potential, feasibility, and impact of these models (providing a detailed representation of P2P transactions). The process and methods used are described in Fig.~\ref{fig:methods}, which comprises: modeling of generation and load patterns for different buses in an electric network; a P2P trade generation process that includes a matching process utilizing a double auction mechanism for buying and selling surplus electricity (assuming the day-ahead market, thus focusing on short-term market integration and scheduling rather than long-term planning), the implementation of proposed P2P trades within a clearance and settlement model based on the outcomes of the matching process and network state for each time frame and lastly, the analyses of results.

\subsection{P2P Trade Generation}
In the absence of established bilateral transactions (P2P), a double auction mechanism is used to generate trades for demonstration purposes. This phase involves simulating generation and load patterns across the electric network's buses to capture the dynamics of electricity demand and supply accurately. Using these simulated results, a double auction mechanism matches buy and sell orders, facilitating P2P trade generation. 

\subsection{Matching process}
The double auction mechanism is used to incorporate real-world constraints on both the demand and supply sides. Let $D$ be the set of buy orders, $S$ the set of sell orders of buses $i \in N$, in a P2P market, whereas:
\begin{itemize}
  \item Buy orders: $D = \{ (d_{i_1} , p^d_{i_1}, q^d_{i_1}), \ldots, (d_{i_m}, p^d_{i_m}, q^d_{i_m}) \mid i \in N, m \in \{1, 2, ..., M_i \}$ where each tuple represents a buy order with identifier $d_i$, a willing-to-pay price $p^d_i$, and a desired quantity $q^d_i$. 
  \item Sell orders: $S = \{ (s_{i_1}, p^s_{i_1}, q^s_{i_1}), \ldots, (s_{i_k}, p^s_{i_k}, q^s_{i_k}) \mid i \in N, k \in \{1, 2, ..., K_i \}$ where each tuple represents a sell order with identifier $s_i$, an asking price $p^s_i$, and an available quantity $q^s_i$.
\end{itemize}
The quantities on each are adjusted at each bus $i \in N$, ${g}_i$ (generation at bus $i$) $\equiv q^s_i$ and adjusted load at bus $i$, ${\rho}_i \equiv q^d_i$ (or that seller cannot sell more than ${g}_i$ and inversely, the buyer cannot buy more than ${\rho}_i$.
The matching process is described in Algorithm~\ref{alg:Matching_algo}.
\begin{algorithm}
\caption{Market Matching Algorithm}
\label{alg:Matching_algo}
\begin{algorithmic}[1]\small
\State $D \gets \text{set of buy orders } \{(d_{i_1},p^d_{i_1},q^d_{i_1}),\ldots,(d_{i_m},p^d_{i_m},q^d_{i_m})\}$
\State $S \gets \text{set of sell orders } \{(s_{i_1},p^s_{i_1},q^s_{i_1}),\ldots,(s_{i_k},p^s_{i_k},q^s_{i_k})\}$
\State $T \gets \emptyset$ \Comment{initialize matches}
\State Sort $D$ by increasing $p^d$
\State Sort $S$ by decreasing $p^s$
\For{each buy order $(d_i,p^d_i,q^d_i)$ in $D$}
    \If{$q^d_i=0$}
        \State \textbf{continue}
    \EndIf
    \For{each sell order $(s_i,p^s_i,q^s_i)$ in $S$}
        \If{$q^s_i>0$ and $p^d_i\geq p^s_i$}
            \State $q_t \gets \min(q^d_i,q^s_i)$
            \State $t \gets (s_i,d_i,q_t)$ \Comment{form trade tuple}
            \State Add $t$ to $T$
            \State $q^d_i \gets q^d_i-q_t$
            \State $q^s_i \gets q^s_i-q_t$
            \If{$q^d_i=0$}
                \State \textbf{break}
            \EndIf
        \EndIf
    \EndFor
\EndFor
\State $Q_{\text{total}} \gets \sum_{t\in T} q_t$ \Comment{Calculate total matched quantity}
\end{algorithmic}
\end{algorithm}
A trade $t \in T$ is represented as $t = (s_t, d_t, q_t)$ where $s_t$ is the seller (source, matching $s_i$ in sell orders), $d_t$ is the buyer (destination, matching $d_i$ in buy orders), and $q_t$ is the quantity of the trade $t$. This ordering scheme simulates an optimal matching: the highest willing buyer is paired with the lowest willing seller. Each user's bid is capped by their respective load for the given time frame, while each ask is constrained by the available local generation. This ensures that the trading process accurately reflects the physical limitations of the energy system. 
\subsection{Clearance and Settlement Model}
The clearance and settlement process within a Peer-to-Peer (P2P) electricity trading framework presents a challenging problem that involves selecting an optimal subset of trades while considering capacity constraints and network structure. Following the establishment of this initial list of proposed trades, the list is then  -- centrally -- processed further to ensure the feasibility of such operations in the network. This is achieved by applying the proposed model (to maximize the total quantity across all trades adhering to network constraints) described in Section~\ref{sec:formulation}, which takes into account the physical constraints and all proposed trades for each time block.

\section{Mathematical formulation}\label{sec:formulation}
\allowdisplaybreaks

\subsection*{Notations and Assumptions}

Assume that time is divided into discrete time blocks. For each time block, the objective is to optimize the volume of P2P trades in the power network. Without loss of generality, we consider a single time block from now. The model is based on the concept that each bus (node) in the power grid can either generate, consume, or transmit
electrical power. In this context, the terms ``bus" and ``node" are used
interchangeably.

Let $N$ be the set of all buses or nodes in the system, where each bus $ i \in N$ is associated with its own power generation and load. The set of nodes $j$ adjacent (directly connected) to node $i$ is denoted by $\text{Adj}(i) = \{ j \mid j \sim i \}$, thus fully characterizing network connectivity and topology. $T$ is the set of all proposed peer-to-peer (P2P) trades $t$ in the network, within $i \in N$. A trade $t \in T$ is represented as $t = (s_t, d_t, q_t)$ where $s_t$ is the seller (source, matching $s_i$ in sell orders of Algorithm~\ref{alg:Matching_algo}), $d_t$ is the buyer (destination matching $d_i$ in buy orders of Algorithm~\ref{alg:Matching_algo}), and $q_t$ is the quantity of the trade $t$. 

For each trade $t$ associated with a potential P2P transaction, we introduce a continuous decision variable $x_t \in [0, 1]$, indicating the extent to which the trade is executed within the network (in other words, we allow for partial execution of trades).

Whenever a trade $t$ is happening, it implies that the destination bus will consume the quantity $x_t q_t$ and the source bus will generate the same quantity $x_t q_t$, as reciprocal. 
Thus, the model also integrates the impact of P2P trades on the power balance at each bus. The trades are factored into the power balance as follows: 

\begin{itemize}
    \item $\sum_{t \in T} x_t q_t \mathbf{1}_{\{s_t = i\}}$, adjusts the generation (${g}_i$), at bus $i$ for trades where $i$ is the source $(s_t = i)$. 
    \item $\sum_{t \in T} x_t q_t \mathbf{1}_{\{d_t = i\}}$, adjusts the load (${\rho}_i$) for bus $i$ for trades where $i$ is the destination $(d_t = i)$.
\end{itemize}

The decision variables $x_t$ determine the execution level of trades, with $q_t$ indicating trade quantities. Indicator functions $\mathbf{1}_{\{s_t = i\}}$ and $\mathbf{1}_{\{d_t = i\}}$ identify node $i$ as the source or destination in a trade $t \in T$, respectively, modulating the trade's impact on the node based on $x_t$ and $q_t$. These functions ensure that each trade's impact is correctly mapped to the respective buses involved. The inflexible (non-dispatchable) generation at bus $i$, not affected by trades, is represented by ${g}_i^0$ and $\sum_{t \in T} x_t q_t \mathbf{1}_{\{s_t = i\}}$ represent (dispatchable generation) adjustments at bus $i$ due to P2P trades, inversely, the inflexible load at bus $i$ is represented by ${\rho}_i^0$ and $\sum_{t \in T} x_t q_t \mathbf{1}_{\{d_t = i\}}$ represent (flexible load) adjustments at bus $i$ due to P2P trades. These operations affect both generation and load depending on the bus's role (seller or buyer) in the trade. $g_i$ and $\rho_i$ are the actual generation and load at bus $i$, after accounting for trades' impact.

Thus, every proposed trade $t$ will adjust $g_i$ and $\rho_i$ on the respective bus:

\begin{equation}
g_i = g^0_i + \sum_{t \in T} x_t q_t \mathbf{1}_{\{s_t = i\}} \ ,
\end{equation}

\begin{equation}
\rho_i = \rho^0_i + \sum_{t \in T} x_t q_t \mathbf{1}_{\{d_t = i\}} \ .
\end{equation}
We use the Direct Current (DC) approximation due to its simplicity \cite{li_numerical_2022} for the power flow constraints, based on voltage angle differences and line susceptance. It simplifies the analysis by focusing on the active power flows and neglecting the reactive power components and losses, as our primary focus is on feasibility and integration of P2P operations. The power flow between any two buses $i$ and $j$ is proportional to the difference in their voltage angles and the line susceptance, capturing the essential dynamics of power transmission in a linearized form. The voltage angle at bus $i$ is denoted by $\theta_i$. $\Delta\theta_{ij} = \theta_i - \theta_j$ is the difference in voltage angles between buses $i$ and $j$, which determines the power flow $P_{ij}$ between them. Considering, $g_i$ and $\rho_i$, into the the nodal power balance ($P_i$) is given by the net generation $g_i$ and consumption $\rho_i$ at each bus $i$, balanced by the sum of power flows to and from its neighboring (adjacent) buses:
\begin{equation}
P_i = (g_i - \rho_i) = \sum_{j \in \text{Adj}(i)} B_{ij} (\theta_i - \theta_j),  \forall i \in N \setminus \{i_\text{ref}\}. \label{eq:pb}
\end{equation}
Thus, the nodal power balance $P$ at each bus $i$ (Equation ~(\ref{eq:pb}) includes the generation, consumption, and power flows to and from adjacent buses, adjusted for any transactions in the network. The susceptance of the line connecting buses $i$ and $j$, as the inverse of the reactance ($X$), is represented as $B_{ij}$. The maximum (absolute) $|P_{ij}|$ power capacity of line connecting buses $i$ and $j$ is denoted by $\overline{P}_{ij}$. It is assumed that $B_{ij} - B_{ji} = 0$. Considering that steady-state stability limit for power transfer, as the susceptance ($B$), is the reciprocal of reactance ($B = \frac{1}{X}$), it can be rewritten as $P \approx B \sin(\Delta \theta)$. The limit used as maximum capacity is set to $\frac{\pi}{6}$ or $30^\circ$ (as a standard operational safety buffer). Additionally, the model includes a slack bus $P_{i_\text{ref}}$, introduced to balance the power in the system and acts as a buffer used to adjust the difference between total generation and total demand within a power system, which is used to balance the power system by adjusting for differences between generation and demand at each node. It is given by:
\begin{equation}
P_{i_\text{ref}} = \sum_{i \in N\backslash\{i_{\text{ref}}\}} (g_i - \rho_i) \ .
\end{equation}

$P_{i_\text{ref}}$ represents the real power component constrained by $P^{\text{ref max}}$, the system's maximum short-circuit power capacity in Mega Volt-Ampere (MVA), assuming a power factor of 1 (purely resistive load). To guarantee system reliability, the absolute magnitude of $ P_{i_\text{ref}}$ (slack bus) is constrained not to surpass its maximum allowed limit:

\begin{equation}
| P_{i_\text{ref}} | \leq P^{\text{ref max}} \ .
\end{equation}

The voltage angle at the reference bus is set to zero, providing a baseline for phase angles: 
\begin{equation}
\theta_{i_\text{ref}} = 0 \ .
\end{equation}
The dual function of bus $P_{i_\text{ref}}$ in a power network can be expressed as follows: absorbing Power (acts as load), when total generation exceeds total load, $P_{i_\text{ref}}$ absorbs excess power, when total load exceeds total generation, $P_{i_\text{ref}}$ supplies deficit power. The sign of $P_{i_\text{ref}}$ determines whether $P_{i_\text{ref}}$ is acting as a load (negative value) or as a generator (positive value).
The trade volume optimization problem is formulated in terms of the decision variables $x_t$, the voltage angles $\theta_i$, and $P_{i_\text{ref}}$. The indicator functions $\mathbf{1}_{\{s_t = i\}}$ and $\mathbf{1}_{\{d_t = i\}}$, are used to determine whether each P2P transactions ($t$) is included and updated on each bus $i$ and are applied to both the source $s_t$ and the destination $d_t$ for each trade $t$ on bus $i$, respectively.

\subsection*{Problem formulation}

With these notations, we can define our model as an optimization problem as follows:

\begin{align}
\underset{x_t, \theta_i, P_{i_\text{ref}}}{\text{maximize}}
& \quad \sum_{t \in T} x_t \cdot q_t \ , \label{eq:objective_2} \\
\text{subject to:} \nonumber \\
& \forall t \in T: 0 \leq x_t \leq 1\ , \label{eq:xt_2} \\
& \forall i \in N \setminus \{i_\text{ref}\}, \forall t \in T: \nonumber \\ 
& g_i = g^0_i + \sum_{t \in T} x_t q_t \mathbf{1}_{\{s_t = i\}}, \label{eq:gen_2} \\
& \rho_i = \rho^0_i + \sum_{t \in T} x_t q_t \mathbf{1}_{\{d_t = i\}}, \label{eq:load_2} \\
& \underline{g}_i \leq {g}_i\leq \overline{g}_i, \label{eq:g_bounds}\\ 
&\underline{\rho}_i \leq {\rho}_i\leq \overline{\rho}_i, \label{eq:rho_bounds} \\
& g_i - \rho_i = \sum_{j \in Adj(i)} B_{ij} (\theta_i - \theta_j), \label{eq:balance_bus_2} \\
& \forall i,j \in N: \nonumber \\
& | B_{ij} (\theta_i - \theta_j) |  \leq \overline{P}_{ij} , \label{eq:capacity_2} \\
& | \theta_i - \theta_j |  \leq \frac{\pi}{6}, \label{eq:angles_2} \\
& \theta_{i_\text{ref}} = 0, \label{eq:anglesref_2} \\
& P_{i_\text{ref}}  = \sum_{i \in N\backslash\{i_{\text{ref}}\}} (g_i - \rho_i), \label{eq:slack_mismatch_2} \\ 
& | P_{i_\text{ref}} | \leq {P^{\text{ref max}}} \ . \label{eq:slack-constraint_2} 
\end{align}

\textbf{Objective Function}
Eqn.~(\ref{eq:objective_2}) represents the objective function that aims to maximize the total quantities of proposed trades, represented by the sum over all trades $t \in T$ of the product of the extent of execution $x_t$ and the respective trade quantity $q_t$.

\textbf{Execution of Trades}
Eqn.~(\ref{eq:xt_2}) is simply our definitional constraint on $x_t$. The continuous decision variables $x_t$ represent the extent to which each trade is executed. These variables directly influence the power balance at the participating buses, and thereby, the overall network's power flow and voltage angles. A higher value of $x_t$ indicates a greater execution of trade $t$, leading to more significant changes in power injections at the corresponding source and destination buses. The execution of trades alters the power injections at specific buses and these changes propagate through the network, affecting voltage angles and power flows. For example, an executed trade increases generation at the source bus ($s_t$) and load at the destination bus ($d_t$), altering the power balance at these buses. This, in turn, affects the voltage angles at these and neighboring buses, thereby influencing the power flows throughout the network.

\textbf{Load and Generation Adjustments Due to Trades}
Eqn.~(\ref{eq:load_2}) combines the initial load and the net effect of trades to provide the total load at each bus. This equation accounts for the adjustments in the load (demand) at each bus due to P2P trades. It adds the net effect of trades to the initial load, effectively updating the demand profile to reflect P2P trading activities. Eqn.~(\ref{eq:gen_2}) combines the initial generation and the net effect of trades to provide the total generation at each bus. This equation accounts for the adjustments in the generation (supply) at each bus due to P2P trades. It adds the net effect of trades to the initial generation, effectively updating the generation profile to reflect P2P trading activities.

\textbf{Generation and load bounds}
Eqn.~(\ref{eq:g_bounds}) denotes the upper $\overline{g}_i$ and lower bounds  $\underline{g}_i$ and Eqn.~(\ref{eq:rho_bounds}), denotes the upper  $\overline{\rho}_i$ and lower bounds  $\underline{\rho}_i$ of ${\rho}_i$.

\textbf{Nodal Power Balance}
Eqn~(\ref{eq:balance_bus_2}) denotes a system of equations where each equation corresponds to a different bus $i$ within the power system, forming a linear system under the DC power flow approximation. Ensures that at each bus $i$, the net power injection $P_i$ ($g_i - \rho_i$) plus the net power exchanged with adjacent buses through P2P trades and transmission lines equals zero. This maintains the principle of power conservation within the network.

\textbf{Line Capacity Constraints}
Eqn.~(\ref{eq:capacity_2}) is the line capacity constraints, ensuring that the power flow between any two connected buses does not exceed the physical limits of the line, maintaining system reliability and preventing overloads.

\textbf{Voltage Angle Differences}
Eqn.~(\ref{eq:angles_2}) ensures that angles remain within a reasonable range. Limits the difference in voltage angles between any two connected buses to prevent excessive phase shifts, which could lead to instability and potentially compromise system operations.

\textbf{Slack Power Mismatch}
Eqn.~(\ref{eq:slack_mismatch_2}) represents the system balance or the total power imported and exported through $P_{i_\text{ref}}$ (or collectively, ``Grid'', as the sum of all activities by other participants). $P_{i_\text{ref}}$ compensates for any mismatches between total generation and total demand, constrained by $P^{\text{ref max}}$, ensuring the overall power balance in the network. 

\section{Implementation}
\subsection{Simulation set up}

To establish a baseline for P2P trade generation, in the residential sector, this work considers NREL's aggregate time-series data \cite{wilson_end-use_2021} for Pennsylvania, ISO regions, and PJM aggregate markets for demand (``End-Use Load Profiles for the U.S. Building Stock''), and Global Solar Atlas \cite{world_bank_group_datainformationmap_nodate} data, assuming a 3kWp solar panel setup, for supply. Hourly demand data undergo K-means clustering to identify typical residential load patterns and solar supply capacity. The networks used for the model are based on the CIGRE \cite{cigre_benchmark_2014} low voltage radial distribution network (44 bus system) and Synthetic Voltage Control LV Networks \cite{lindner_aktuelle_2016} ``Village'' (80 bus system).

\begin{figure}[htbp]
   \vspace{-4mm}
  \begin{subfigure}{0.485\columnwidth}
    \centering
    \includegraphics[width=\linewidth]{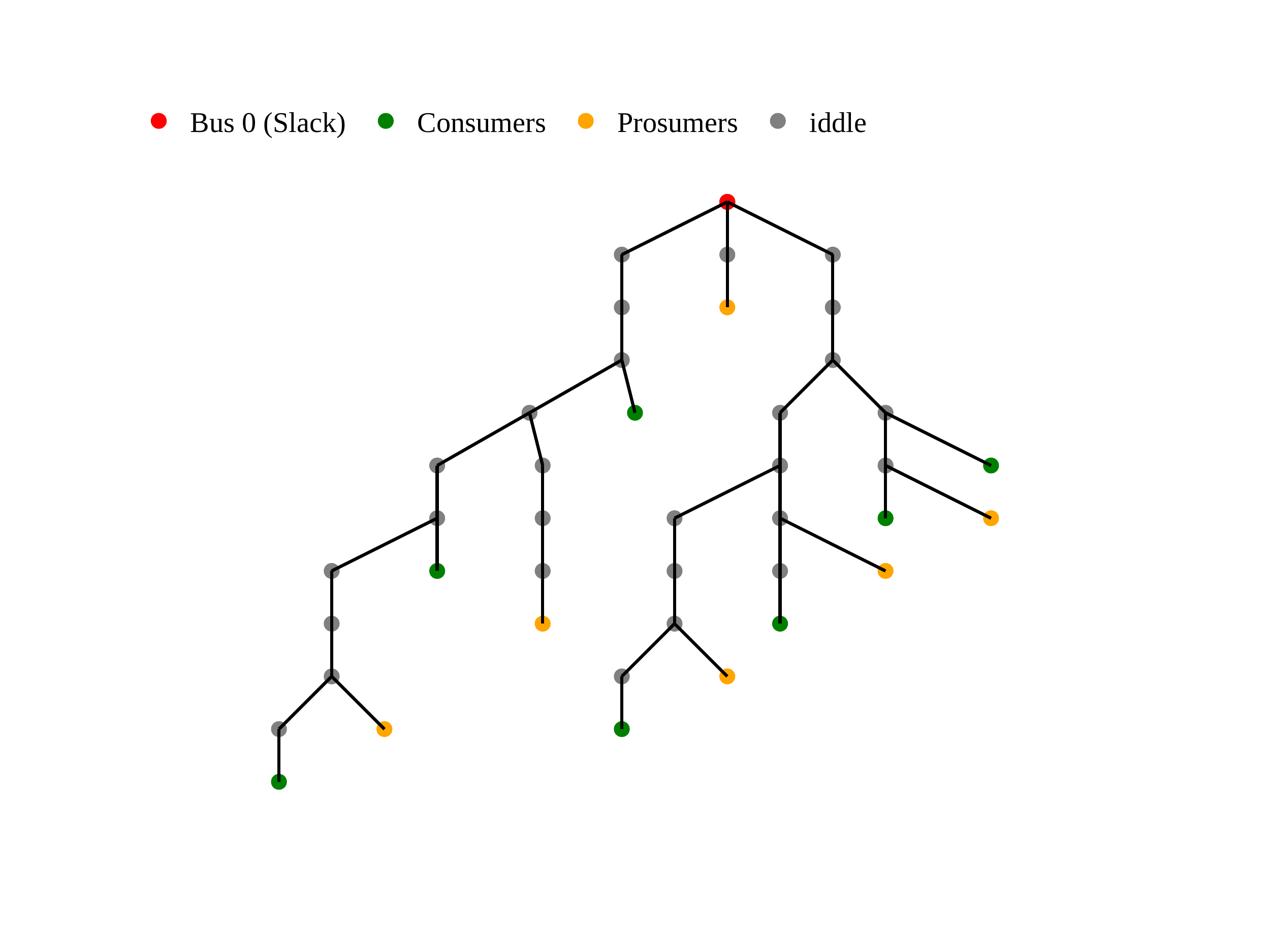}
    \caption{Adapted CIGRE \cite{cigre_benchmark_2014} low voltage radial distribution network (44 bus system) with users' profile}
    \label{fig:cigre}
  \end{subfigure}
  \begin{subfigure}{0.485\columnwidth}
    \centering
    \includegraphics[width=\linewidth]{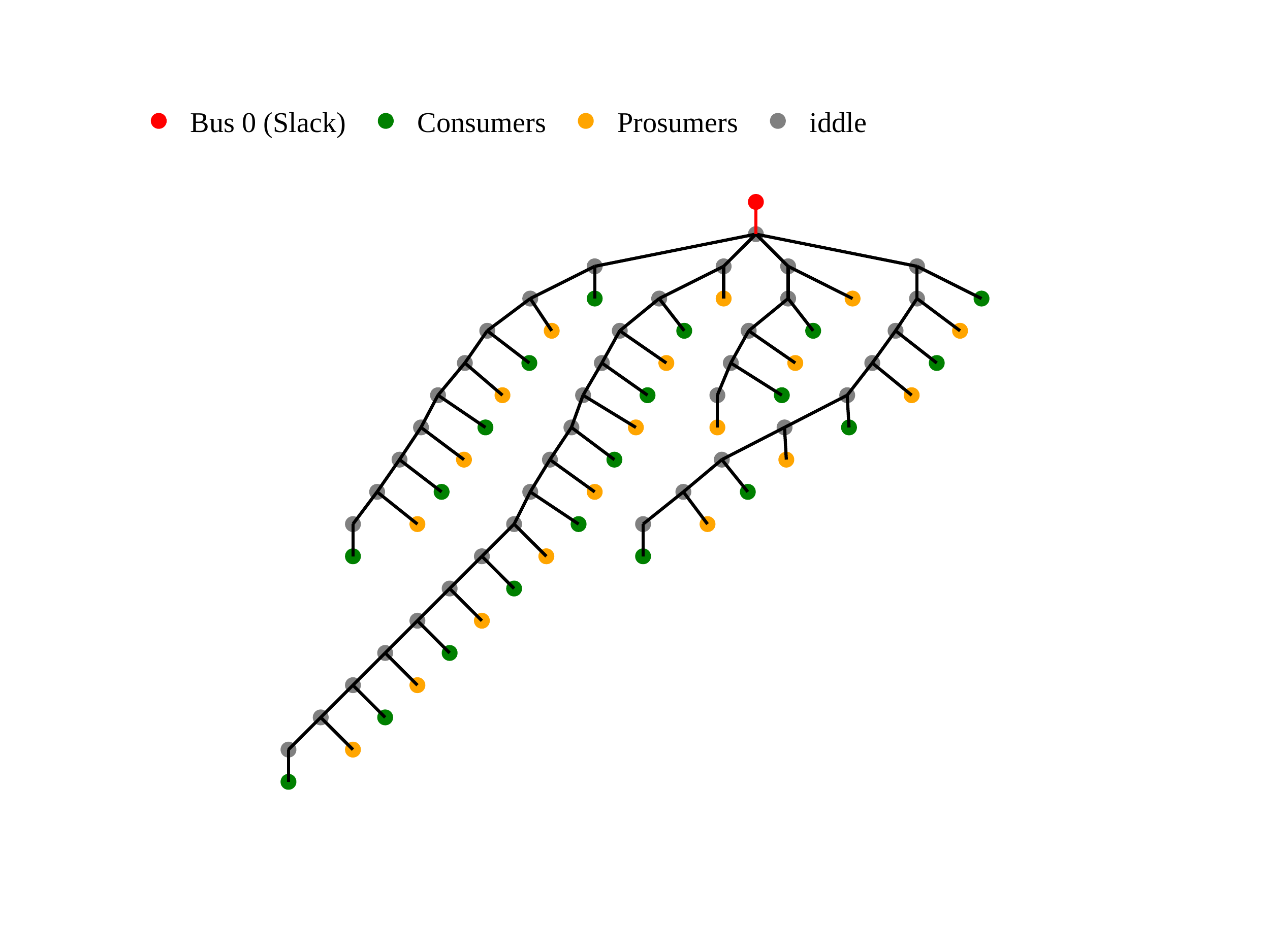}
    \caption{Adapted Synthetic Voltage Control LV Networks \cite{lindner_aktuelle_2016} "Village" (80 bus system) with users' profile}
    \label{fig:net_village}
  \end{subfigure}
    \vspace{-2mm}
\end{figure}

For local supply, 50\% of overall consumers (buses) are associated with local solar generation capacity on the same bus to look at different local generation profiles (see Figures~\ref{fig:cigre} and~\ref{fig:net_village}). In these figures, the orange-colored nodes represent consumers with associated local solar generation, the green nodes represent consumers without any associated generation, and the red corresponds to the $P_{i_\text{ref}}$ (LV/MV connection).
\begin{figure}[htbp]
\vspace{5mm} 
  \begin{subfigure}{0.49\linewidth}
    \centering
    \includegraphics[width=\linewidth]{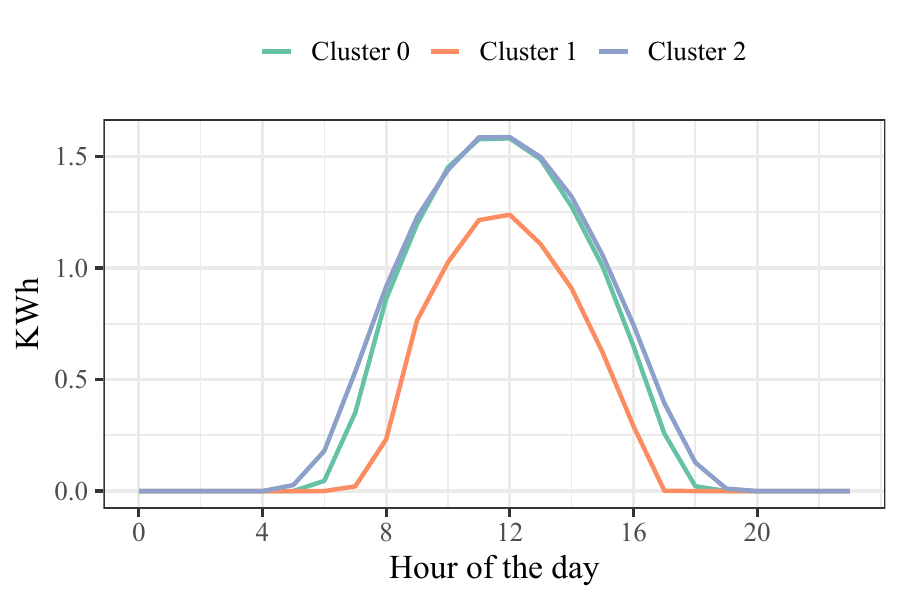}
    \caption{Solar PV}
    \label{fig:solar_clusters} 
  \end{subfigure}
  \begin{subfigure}{0.49\linewidth}
    \centering
    \includegraphics[width=\linewidth]{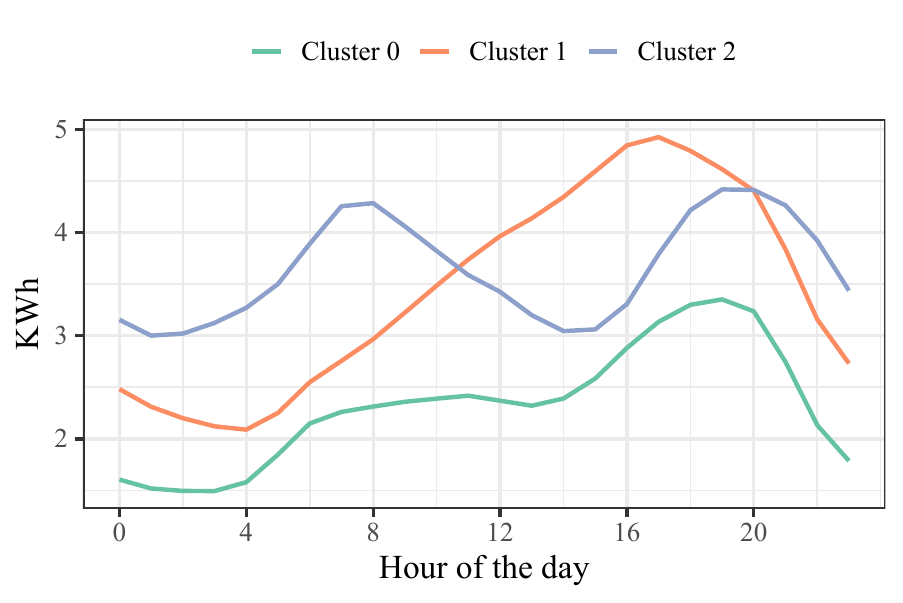}
    \caption{Residential loads} 
    \label{fig:clusters_residential}
  \end{subfigure}
  \caption{Daily clusters}
  \vspace{-5mm}
\end{figure}

In the simulation, the double auction mechanism is used to incorporate real-world constraints on both the demand and supply sides. Each user's bid is capped by their respective load for the given time block, while each ``ask'' is constrained by the available local solar generation (as described in Section~\ref{sec:methods}). It was considered that the reserve prices for each side were uniform and constant across each time block. Following the establishment of this initial set of proposed transactions is then processed further by applying the model (to maximize the total quantity across all trades adhering to network constraints) described earlier in Section~\ref{sec:formulation}, which takes into account the physical constraints. This ensures that the resulting solution (accepted trades) is not only efficient within each side of the P2P exchange but also physically feasible. In the case of full P2P, this process is rather simplified as there is already matching offers (or a set of bilateral contracts, P2P trades $t \in T$). The binary decision variables are relaxed to take on real values between 0 and 1, and $x_t$ acts as the scaling factor, reformulating as a Linear Program (LP). As an LP the problem can be solved by Simplex method or Interior-Point-Method (IPM), where the last was used with Gurobi solver. Pandapower library \cite{thurner_pandapoweropen-source_2018} was employed for the manipulation of power system data.\footnote{The code is available at \url{https://github.com/d-vf/P2PEnergyTrading}} Lastly, the extreme case was set at $10^2$ of baseline, to evaluate the ability to accommodate extreme values and their (expected) behavior.

\section{Results}
The analysis is presented in two distinct tables, each representing key metrics from a 24-hour block analysis: one for a 44-bus system (Table~\ref{tab:results_44}) and another for an 80-bus system (Table~\ref{tab:results_80}).``Trade Qty. (kWh)'' is the result of the optimization process, representing the maximum sum of all $x_t \cdot q_t$, where  $x_t$ is the decision variable for trade acceptance, and $q_t$ is the trade quantity for each seller/generator $i$ in trade $t$ for each hour, ``slack power (kWh)'' indicates the power managed by the $P_{i_\text{ref}}$, functioning as an absorber or generator to maintain system balance in an open system. ``Loads (kWh)'' shows the total power demand across all buses per hour. ``Count trades'' enumerates the total number of permissible trades and ``\% Load Fulfilled w/ P2P'' reflect the percentage of electricity fulfilled trough P2P transactions.
\begin{table}[ht!]
\centering
\vspace{1mm} 
\caption{ 44 bus system -- 24 hour period}
\begin{adjustbox}{width=0.99\columnwidth,center}
  \begin{tabular}{c p{1.5cm} p{1.5cm} p{1.5cm} p{1.5cm} p{1.5cm}}
  \toprule
  \textbf{Hour} & \textbf{Trade Qt (kWh)} & \textbf{Count Trades (t=6)} & \textbf{Slack Power (kWh)} & \textbf{Loads (kWh)} & \textbf{\% Load Fulfilled w/ P2P} \\
  \midrule
  0--5 & -- & -- & -- & -- & -- \\
  6 & $< 0.01$ & 6 & -12.898 & 12.898 & $< 1\%$ \\
  7 & $0.126$ & 6 & -13.443 & 13.569 & 0.928\% \\
  8 & $1.392$ & 6 & -12.486 & 13.878 & 10.03\% \\
  9 & $4.596$ & 6 & -9.561 & 14.157 & 32.464\% \\
  10 & $6.156$ & 6 & -8.185 & 14.341 & 42.923\% \\
  11 & $7.290$ & 6 & -7.222 & 14.512 & 50.231\% \\
  12 & $7.434$ & 6 & -6.787 & 14.221 & 52.274\% \\
  13 & $6.642$ & 6 & -7.286 & 13.928 & 47.685\% \\
  14 & $5.448$ & 6 & -8.898 & 14.346 & 37.974\% \\
  15 & $3.732$ & 6 & -11.773 & 15.505 & 24.069\% \\
  16 & $1.746$ & 6 & -15.542 & 17.288 &10.098\% \\
  17 & $0.060$ & 6 & -18.801 & 19.792 & 0.319\% \\
  18-23  & -- & -- & -- & -- & -- \\
  \bottomrule
  \end{tabular}%
\end{adjustbox}
\label{tab:results_44}
\end{table}
\begin{figure}
  \centering
  \includegraphics[width=0.9\columnwidth]{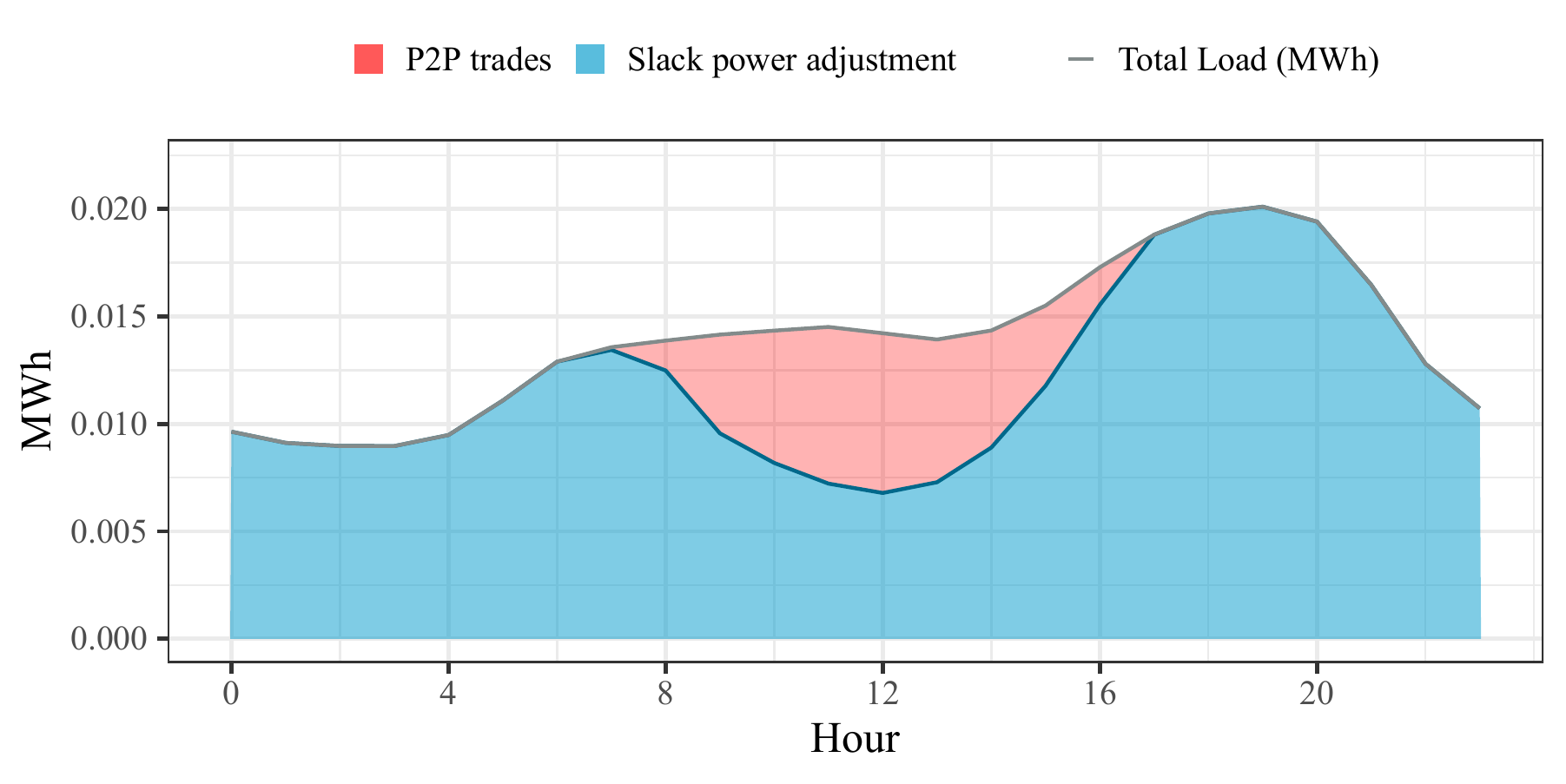} \caption{Total P2P traded and slack adjustments per hour -- 44 bus system}
  \label{fig:stacked}
\vspace{-6mm}
\end{figure}
In this simulated scenario, the quantity available for trade is determined by the output of the generators (or sellers), specifically those using solar generation (P2P Fig. \ref{fig:stacked}), where there is higher integration of renewable electricity at hours with most solar output (e.g. more than 50\% at 11-hour block for both scenarios).

\begin{table}[ht!]
\centering
\vspace{1mm}
\caption{ 80 bus system -- 24 hour Period}
  \begin{adjustbox}{width=0.99\columnwidth,center}
    \begin{tabular}{c p{1.5cm} p{1.5cm} p{1.5cm} p{1.5cm} p{1.5cm}}
      \toprule
      \textbf{Hour} & \textbf{Trade Qt (kWh)} & \textbf{Count Trades (t=20)} & \textbf{Slack Power (kWh)} & \textbf{Loads (kWh)} & \textbf{\% Load Fulfilled w/ P2P} \\
      \midrule
      0--5 & -- & -- & -- & -- & -- \\
      6 & $< 0.01$  & 20 & -42.994 & 42.994 & $< 1 \%$ \\
      7 & $0.420$ & 20 & -44.791 & 45.232 & 0.92 \% \\
      8 & $4.640$ & 20 & -41.388 & 46.260 & 10.00 \% \\
      9 & $15.320$ & 20 & -31.104 & 47.190 & 32.99 \% \\
      10 & $20.520$ & 20 & -26.259 & 47.805 & 43.86 \% \\
      11 & $24.300$ & 20 & -22.860 & 48.375 & 51.52 \% \\
      12 & $24.780$ & 20 & -21.384 & 47.403 & 48.84  \% \\
      13 & $22.140$ & 20 & -23.182 & 46.429 & 38.70 \% \\
      14 & $18.160$ & 20 & -28.753 & 47.821 & 24.36 \% \\
      15 & $12.440$ & 20 & -38.622 & 51.684 & 24.62 \% \\
      16 & $5.820$ & 20 & -51.518 & 57.629 & 10.15 \% \\
      17 & $0.020$ & 20 & -62.670 & 62.691 & 0.310 \% \\
      18--23 & -- & -- & -- & -- & -- \\
      \bottomrule
    \end{tabular}%
  \end{adjustbox}
  \label{tab:results_80}
\vspace{-1mm}
\end{table}

Table~\ref{tab:extreme_44} shows the extreme case -- aiming to test the robustness of the proposed model -- set at $10^2$ proposed P2P of the baseline, for the 11-hour block, to evaluate the ability to accommodate extreme proposed trade volume and their behavior. In this specific scenario, the reduction in $x_t$ identified as a relaxed binary variable signifies the system's modifications -- decrease in the quantity (or partial amounts) of transactions permitted -- to uphold network constraints. $P_{i_\text{ref}}$ denotes the electricity Imported or Exported (Slack I/E) for the same time block.
\begin{table}[ht!]
  \centering
  \vspace{-1mm} 
  \caption{Extreme P2P Trades -- 44 bus system}  
  \setlength{\belowcaptionskip}{-20pt} 
  \setlength{\belowdisplayskip}{-20pt} 
  \begin{adjustbox}{width=0.8\columnwidth,center}
    \begin{tabular}{cccc}
      \toprule
      \textbf{Hour} & \textbf{Seller (Bus)} & \textbf{$x_t$} & \textbf{Power (kWh)} \\
      \midrule
      11 & Seller 0 (Bus 16) & 0.534 & 64.946 \\
      11 & Seller 1 (Bus 18) & 0.524 & 63.742 \\
      11 & Seller 2 (Bus 22) & 0.009 &  1.204 \\
      11 & Seller 3 (Bus 36) & 0.344 & 41.887 \\
      11 & Seller 4 (Bus 40) & 0.344 & 41.887 \\
      11 & Seller 5 (Bus 42) & 0.344 & 41.887 \\
      \midrule
      \multicolumn{2}{c}{\textbf{Total Load}} & \multicolumn{2}{c}{262.776} \\
      \midrule
      \multicolumn{2}{c}{\textbf{Total P2P Generation}} & \multicolumn{2}{c}{255.552} \\
      \midrule
      \multicolumn{2}{c}{\textbf{$P_{i_\text{ref}}$ (Slack I/E)}} & \multicolumn{2}{c}{-07.223} \\
      \bottomrule
    \end{tabular}%
  \end{adjustbox}
  \label{tab:extreme_44} 
\vspace{-1mm}
\end{table}

\section{Conclusion}
We developed and formalized a model addressing the clearance and settlement problem in peer-to-peer (P2P) trading within low-voltage networks, with multiple users and different time windows which optimizes electricity P2P trades and ensures reliable and transparent transactions. The presented model validates a dispatch ``plan'' or ``schedule'' of feasible P2P trades that the network can handle without violating any operational constraints, similar to the process in (wholesale) market operations, but doing it preemptively. This forward-looking approach not only maximizes unused local capacity (from P2P) but also ensures that the electrical system remains stable and within its operational boundaries when actual dispatched traded P2P occurs. As shown in the carried simulations, it allows the integration of unused close capacity and excess renewable local generation (e.g., solar generation). The model is designed to secure the fulfillment of dispatchable loads, effectively merging these resources with existing market mechanisms, in an open manner, with multiple users. Conversely, in an extreme case scenario, the model is tested under conditions where there is a significantly high quantity available for trade, allowing for an exploration of the system's behavior under unusually high trading volumes. The proposed model contributes to the development and implementation of new market structures, such as P2P, by creating clearance and settlement mechanisms in low-voltage networks, critical elements to allow the existence of active participation of ``prosumers'' while keeping the resulting power flow within acceptable bounds. This model can also be adapted to energy communities by implementing an adjusted sharing coefficient of allocation or simply through aggregation operations targeting the off-taker, instead of direct P2P allocation. There is a variety of sources from which consumers can acquire electricity -- ranging from P2P to local energy communities to traditional retail and grid suppliers -- that can be integrated within network operations.

\bibliography{references_zotoro} 

\begin{thebibliography}{10}
\bibitem{world_bank_group_datainformationmap_nodate}
W.~B. Group, ``[{Data}/information/map] obtained from the “{Global} {Solar} {Atlas} 2.0, a free, web-based application is developed and operated by the company {Solargis} s.r.o. on behalf of the {World} {Bank} {Group}, utilizing {Solargis} data, with funding provided by the {Energy} {Sector} {Management} {Assistance} {Program} ({ESMAP})..''
\end{thebibliography}

\bibliographystyle{ieeetr} 

\end{document}